\documentclass{CSML}

\pdfoutput=1

\usepackage{lastpage}

\def\dOi{13(3:30)2017}
\lmcsheading%
{\dOi}
{1--\pageref{LastPage}}
{}
{}
{Jan.~13, 2014}
{Sep.~26, 2017}
{}


\usepackage[hidelinks]{hyperref}

\RequirePackage[latin1]{inputenc}
\RequirePackage{url}
\newcommand{
 \ifthenels}[3]{[\textbf{if}\; #1\; \textbf{then}\; #2\; \textbf{else}\; #3]}

\newcommand{\C}[1]{\mathcal{#1}}

\newcommand{\ex}{\exists}
 
\newcommand{\exnc}{\exists^{\mathrm{u}}}

\newcommand{\pair}[2]{\langle #1, #2 \rangle}
\newcommand{\grec}{\C{F}}
\newcommand{\typeN}{\mathbf{N}}
\newcommand{\nullterm}{\varepsilon}
\newcommand{\defeq}{:=}
\newcommand{\defequiv}{:=}
\newcommand{\inquotes}[1]{``#1''}
\newcommand{\Prog}{\mathrm{Prog}}

\theoremstyle{plain}

\begin{document}

\title[A bound for Dickson's lemma]{A bound for Dickson's lemma}

\author[J.~Berger]{Josef Berger}
\address{Mathematisches Institut der LMU
  M\"unchen, Theresienstra\ss{}e 39, 80333 M\"unchen}
\email{\{jberger,schwicht\}@math.lmu.de}  

\author[H.~Schwichtenberg]{Helmut Schwichtenberg}
\address{\vspace{-18 pt}}

\begin{abstract}
  We consider a special case of Dickson's lemma: for any two functions
 $f,g$ on the natural numbers there are two numbers $i<j$ such that
 both $f$ and $g$ weakly increase on them, i.e., $f_i\le f_j$ and
 $g_i \le g_j$.  By a combinatorial argument (due to the first
 author) a simple bound for such $i,j$ is constructed.  The
 combinatorics is based on the finite pigeon hole principle and
 results in a descent lemma.  From the descent lemma one can prove
 Dickson's lemma, then guess what the bound might be, and verify it
 by an appropriate proof.  We also extract (via realizability) a
 bound from (a formalization of) our proof of the descent lemma.
\end{abstract}

\maketitle

\section{Introduction}\label{S:Introduction}
Consider the following special case of Dickson's lemma: for any two
functions $f,g$ on the natural numbers there are two numbers $i<j$
such that both $f$ and $g$ weakly increase on them, i.e., $f_i\le f_j$
and $g_i \le g_j$.  By a combinatorial argument (due to the first
author) a simple bound for such $i,j$ is constructed.  The
combinatorics is based on the finite pigeon hole principle and results
in a certain descent lemma.  From the descent lemma one can prove
Dickson's lemma, then directly guess what the bound might be, and
finally verify it by an appropriate proof.  We also extract (via
realizability) a bound from (a formalization of) our proof of the
descent lemma.

In its usual formulation, Dickson's lemma (for fixed functions) is a
$\Sigma^0_1$-formula.  In contrast, we shall prove a quantifier-free
statement which implies Dickson's lemma in its usual form, but not
vice versa.  Our proof can be carried out in the formal system of
Elementary Analysis \cite[p.144]{TroelstravanDalen88}, a conservative
extension of Heyting arithmetic with variables and quantifiers for
number-theoretic functions.  In fact, we don't make use of the axiom
of choice at all.  Furthermore, we can restrict induction to
quantifier-free formulas.

Dickson's lemma has many applications.  For instance, it is used to
prove termination of Buchberger's algorithm for computing Gr\"obner
bases \cite{Buchberger70}, and to prove Hilbert's basis theorem
\cite{Simpson88}.

There are many other proofs of Dickson's lemma in the literature, both
with and without usage of non-constructive (or \inquotes{classical})
arguments.  The original proof of Dickson \cite{Dickson13} and the
particularly nice one by Nash-Williams \cite{Nash-Williams63} (using
minimal bad sequences) are non-constructive, and hence do not
immediately provide a bound.  But it is well known that by using some
logical machinery one can still read off bounds, using either
G\"odel's \cite{Goedel58} Dialectica translation as in Hertz
\cite{Hertz04} or Friedman's \cite{Friedman78} $A$-translation as in
\cite{BergerBuchholzSchwichtenberg02}.  However, these bounds -- even
for the case of just two functions considered here -- heavily use
higher type (primitive recursive) functionals and are less perspicious
than the one obtained below.

The first constructive proof of Dickson's lemma has been given by
Sch\"utte and Simpson \cite{SchuetteSimpson85,Simpson88}, using
ordinal numbers and transfinite induction up to $\epsilon_0$.  Similar
methods have been used by Sustik \cite{Sustik03} and Mart\'{i}n-Mateos
et al.\ \cite{MMRRAH11}.  Since initial segments of transfinite
induction are used, these proofs when written in arithmetical systems
require ordinary induction on quantified formulas.  A different
constructive proof has been given by Veldman \cite{Veldman04}.  It
uses
dependent choice for $\Sigma_1$-formulas (with parameters), and
induction on $\Pi_2$-formulas.  This proof also provides the basis of
Fridlender's \cite{Fridlender97} formalization in Agda.  The
computational content of these proofs has not been studied; the bound
involved will be very different from the present one.

\section{A combinatorial proof of Dickson's lemma}
\label{S:Dickson}
We start with a finite pigeonhole principle, in two disjunctive forms.
The (rather trivial) proofs are carried out because they have
computational content which will influence the term extracted from a
formalization of our proofs in Section~\ref{S:Extract}.

\begin{lem}[FPHDisj]
 \label{L:FPHDisj}
 $\forall_{m,f}(\ex_{i<j \le m} f_i=f_j \lor \ex_{j \le m} m \le f_j)$.
\end{lem}

\begin{proof}
 By induction on $m$.  For $m=0$ the second alternative holds.  In
 case $m+1$ let $f_j$ be maximal among $f_0, \dots, f_{m+1}$.  If
 $m+1 \le f_j$ we are done.  Else we have $f_j \le m$.  Now we apply
 the induction hypothesis to $f' \defeq f_0, \dots, f_{j-1}, f_{j+1},
 \dots, f_{m+1}$.  If two of them are equal we are done.  Else $m \le
 f_k$ for some $k \ne j$ and hence $f_j \le f_k$.  If $f_j=f_k$ we
 are done.  Else we have $f_j<f_k$, contradicting the choice of $j$.
\end{proof}

Note that quantifier-free induction suffices here, since we only prove
a property of finite lists of natural numbers.

In the key lemma~\ref{L:Key} below we will need a somewhat stronger
disjunctive version of the pigeonhole principle.  To this end we need
an injective coding $\pair {n} {m}$ of natural numbers which is
\inquotes{square-filling}, i.e.\ with the property
\begin{equation}
 \label{E:CodeSqFill}
 k^2 \le \pair {n} {m} \to k \le n \lor k \le m.
\end{equation}
This can be achieved by
\begin{alignat*}{5}
 &\dots \\
 &12 &\quad& 13 &\quad& 14  &\quad& 15 &\quad& \dots\\
 &6  &\quad&  7 &\quad& 8 &\quad& 11 &\quad&  \dots\\
 &2  &\quad&  3 &\quad& 5 &\quad& 10 &\quad& \dots\\
 &0  &\quad&  1 &\quad& 4 &\quad& 9 &\quad& \dots
\end{alignat*}
or explicitely
\begin{equation*}
 \pair {n} {m} \defeq
 \begin{cases}
   n^2+m &\hbox{if $m<n$},
   \\
   m^2+m+n &\hbox{otherwise}.
 \end{cases}
\end{equation*}

\begin{lem}[FPHDisj2]
 \label{L:FPHDisj2}
 \begin{equation*}
   \forall_{f,g,k}(
   \ex_{i<j \le k^2}(f_i=f_j \land g_i=g_j) \lor
   \ex_{j \le k^2} k \le f_j \lor
   \ex_{j \le k^2} k \le g_j).
 \end{equation*}
\end{lem}

\begin{proof}
 Fix $f,g,k$.  Use Lemma~\ref{L:FPHDisj} with $s_i \defeq \pair {f_i}
 {g_i}$ and $m \defeq k^2$.  In the first case from $s_i=s_j$ we
 obtain $f_i=f_j$ and $g_i=g_j$ by the injectivity of the coding.  In
 the second case we have some $j \le k^2$ with $k^2 \le s_j$.  From
 the square-filling property \eqref{E:CodeSqFill} of the coding we
 obtain $k \le f_j$ or $k \le g_j$.
\end{proof}

As an immediate consequence we have

\begin{lem}[Key]
 \label{L:Key}
 \begin{align*}
   \forall_{f,g,n,k}(
   &\ex_{n<i<j \le n+k^2+1}(f_i=f_j \land g_i=g_j) \lor {}
   \\
   &\ex_{n<j \le n+k^2+1} k \le f_j \lor
   \ex_{n<j \le n+k^2+1} k \le g_j).
 \end{align*}
\end{lem}

\begin{proof}
 Use Lemma~\ref{L:FPHDisj2} for $\lambda_i f_{n+1+i}$, $\lambda_i
 g_{n+1+i}$ and $k$.
\end{proof}

Now we introduce some notation.  $\mathrm{Mini}(f,n)$ is the first
argument where $f$ is minimal on $\{0, \dots, n\}$:
\begin{align*}
 &\mathrm{Mini}(f,0) \defeq 0,
 \\
 &\mathrm{Mini}(f,n+1) \defeq
 \begin{cases}
   \mathrm{Mini}(f,n)
   &\hbox{if $f_{\mathrm{Mini}(f,n)} \le f_{n+1}$},
   \\
   n+1
   &\hbox{otherwise}.
 \end{cases}
\end{align*}
We define functions $\Psi, \Phi, I$ and a formula $D$ with arguments
$f,g,n$.  For readability $f,g$ are omitted.
\begin{equation}
 \label{E:PsiPhiID}
 \begin{split}
   &\Psi(n) \defeq \max \{ f_{\mathrm{Mini}(g,n)}, g_{\mathrm{Mini}(f,n)} \},
   \\
   &\Phi(n) \defeq f_{\mathrm{Mini}(f,n)}+g_{\mathrm{Mini}(g,n)},
   \\
   &I(n) \defeq n+\Psi(n)^2 +1,
   \\
   &D(n) \defeq \ex_{i<j\le n}(f_i\le f_j \land g_i \le g_j).
 \end{split}
\end{equation}
$D(n)$ expresses that $n$ is a bound for Dickson's lemma.

The next lemma states a crucial property of the function $I$: either
$I(n)$ already is a bound for Dickson's lemma, or else $\Phi$
decreases properly when going from $n$ to $I(n)$.  Since this cannot
happen infinitely often, iteration of $I$ will finally give us the
desired bound.

\begin{lem}[Descent]
 \label{L:Descent}
 $D(I(n)) \lor \Phi(I(n))<\Phi(n)$.
\end{lem}

\begin{proof}
 Use Lemma~\ref{L:Key} with $f$, $g$, $n$ and $\Psi(n)$.  In the first
 case we have $D(I(n))$.  In the second case we have $n<j \le I(n)$
 with $\Psi(n) \le f_j$; the third case is symmetric.  Let $i \defeq
 \mathrm{Mini}(g,n)$.  Then $f_i \le \Psi(n)$.  In case $g_i \le g_j$
 we have $D(I(n))$ and are done.  Therefore assume $g_j<g_i$.  We
 show (i) $\Phi(I(n)) \le \Phi(j)$ and (ii) $\Phi(j)<\Phi(n)$.  From
 $j \le I(n)$ we obtain (i).  For (ii) we show
 $f_{\mathrm{Mini}(f,j)}+g_{\mathrm{Mini}(g,j)} <
 f_{\mathrm{Mini}(f,n)}+g_i$.  Now $n<j$ implies
 $f_{\mathrm{Mini}(f,j)} \le f_{\mathrm{Mini}(f,n)}$, and
 $g_{\mathrm{Mini}(g,j)} \le g_j<g_i$.
\end{proof}

From Lemma~\ref{L:Descent} we construct a bound for Dickson's lemma.
Let
\begin{equation*}
 I^0(n) \defeq n,\quad I^{m+1}(n) \defeq I(I^m(n)).
\end{equation*}

\begin{lem}
 \label{L:DescentCor}
 $D(I^n(0)) \lor \Phi(I^n(0))+n \le \Phi(0)$.
\end{lem}

\begin{proof}
 Induction on $n$.  Step $n \mapsto n+1$.  Applying Lemma~\ref{L:Descent}
 to $I^n(0)$ gives $D(I^{n+1}(0)) \lor
 \Phi(I^{n+1}(0))<\Phi(I^n(0))$.  In the second case we have
 \begin{equation*}
   \Phi(I^{n+1}(0))+n+1 < \Phi(I^n(0))+n+1
   \le \Phi(0)+1
 \end{equation*}
 The latter inequality follows from the induction hypothesis, since
 $D(I^n(0))$ implies $D(I^{n+1}(0))$.
\end{proof}

\begin{prop}
 \label{P:BoundDickson}
 $D(I^{f_0+g_0+1}(0))$.
\end{prop}

\begin{proof}
 Apply Lemma~\ref{L:DescentCor} to $\Phi(0)+1$.
\end{proof}

This bound is far from optimal: already for
\begin{equation*}
 f_n \defeq \begin{cases} 1 &\hbox{if $n=0$},\\ 0 &\hbox{else} \end{cases}
 \qquad
 g_n \defeq 0
\end{equation*}
with optimal bound $2$ we have
\begin{equation*}
 I^{f_0+g_0+1}(0) = I^2(0) =I(I(0)) > I(0) = \Psi(0)^2 +1 = 2.
\end{equation*}

Can we extend this proof to show Dickson's lemma for finitely many
functions?  For instance for three functions a corresponding version
of the key lemma holds:
\begin{align*}
 \forall_{f,g,h,n,k}(
 &\ex_{n<i<j \le n+k^4+1}(f_i=f_j \land g_i=g_j \land h_i=h_j) \lor {}
 \\
 &
 \ex_{n<j \le n+k^4+1} k \le f_j \lor
 \ex_{n<j \le n+k^4+1} k \le g_j \lor
 \ex_{n<j \le n+k^4+1} k \le h_j)
\end{align*}
(Proof.  Apply the original key lemma to $\pair f g, h, n$ and $k^2$).
We can also define a measure function $\Phi(n) \defeq
f_{\mathrm{Mini}(f,n)} + g_{\mathrm{Mini}(g,n)} + h_{\mathrm{Mini}(h,n)}$.
A natural candidate for $\Psi$ is
\begin{equation*}
 \Psi(n) \defeq \max\{
 f_{\mathrm{Mini}(g,n)},
 f_{\mathrm{Mini}(h,n)},
 g_{\mathrm{Mini}(f,n)},
 g_{\mathrm{Mini}(h,n)},
 h_{\mathrm{Mini}(f,n)},
 h_{\mathrm{Mini}(g,n)} \}
\end{equation*}
and a natural candidate for $I$ is $I(n) \defeq n + \Psi(n)^4 +1$.
But the corresponding version of the descent lemma is false: let $n
\defeq 2$ and
\begin{align*}
 f \defeq (0,1,1,1,0,f_5,\dots),
 \\
 g \defeq (1,0,1,0,1,g_5,\dots),
 \\
 h \defeq (1,1,0,0,0,h_5,\dots).
\end{align*}
Then $\Phi(n)=0$, $\Psi(n)=1$, $I(n)=4$, and we neither have $D(I(n))$
nor $\Phi(I(n))<\Phi(n)$. -- However, it may well be that a more
refined form of the present approach works.  We leave this for future
research.

\section{Extracting computational content}
\label{S:Extract}
In the following, we demonstrate how a bound for Dickson's lemma can
be extracted from a proof of the existence of such a bound.  The proof
we will use is essentially the one presented in
Section~\ref{S:Dickson}, i.e., it is based on the descent
lemma~\ref{L:Descent}.  We will then apply the realizability
interpretation to obtain the bound.  In fact, the bound will be
\emph{machine} extracted from a formalization of the existence proof.

In more detail, we shall use that $I$ is increasing (i.e., $n<I(n)$)
and that from $D(n)$ and $n<m$ we can infer $D(m)$.  Then we prove the
existence of a bound by general induction with measure $\Phi$.

\subsection{General induction and recursion}
\label{SS:GenRec}
We first explain general induction w.r.t.\ a measure, and the
corresponding definition principle of general recursion.

General induction allows recurrence to all points \inquotes{strictly
 below} the present one.  In applications it is best to make the
necessary comparisons w.r.t.\ a measure function $\mu$; for simplicity
we restrict ourselves to the case where $\mu$ has values in the
natural numbers, and the ordering we refer to is the standard
$<$-relation.  The principle of general induction then is
\begin{equation*}
 \forall_{\mu, x} ( \Prog^{\mu}_x P x \to P x ),
\end{equation*}
where $\Prog^{\mu}_x P x$ expresses \inquotes{progressiveness}%
\index{progressive} w.r.t.\ $\mu$ and $<$, i.e.,
\begin{equation*}
 \Prog^{\mu}_x P x \defequiv
 \forall_x (
 \forall_y(\mu y < \mu x \to P y) \to P x ).
\end{equation*}
It is easy to see that in our special case of the $<$-relation we can
prove  general induction from structural induction.  However, it will be
convenient to use general induction as a primitive axiom, for then the
more efficient general recursion constant $\grec$ will be extracted.
It is defined by
\begin{equation*}
 \grec \mu x G = G x (
 \lambda_y \ifthenels{\mu y<\mu x}{\grec \mu y G}{\nullterm}),
\end{equation*}
where $\nullterm$ denotes a canonical inhabitant of the range.  It is
easy to prove that $\grec$ is definable from an appropriate structural
recursion operator.

\subsection{Non-computational quantifiers}
\label{SS:NcQuant}
We now use general induction in our constructive proof of Dickson's
lemma.  However, we have to be careful with the precise formulation of
what we want to prove.  We are not interested in the pair $i,j$ of
numbers where both $f$ and $g$ increase, but only in a bound telling
us when at the latest this must have happened.  Therefore the
existential quantifiers $\ex_{i,j}$ must be made \inquotes{uniform}
(i.e., non-computational); it will be disregarded in the realizability
interpretation.  Such non-computational quantifiers have first been
introduced in \cite{Berger93a,Berger05}; in
\cite{SchwichtenbergWainer12} this concept is extended to all
connectives and discussed in detail.  Let
\begin{equation*}
 D'(n) \defeq \exnc_{i<j\le n}(f_i\le f_j \land g_i \le g_j).
\end{equation*}
Using this non-computational form of $D(n)$ we modify Lemma~\ref{L:Descent}
to
\begin{lem}[Descent$^{\mathrm{nc}}$]
 \label{L:DescentModif}
 $D'(I(n)) \lor \Phi(I(n))<\Phi(n)$.
\end{lem}
Note that the computational content of a proof of this lemma is that
of a functional mapping two unary functions and a number into a
boolean.  From Lemma~\ref{L:DescentModif} we obtain as before a
modification of Proposition~\ref{P:BoundDickson} to
\begin{prop}[Bound for Dickson's lemma]
 \label{P:BoundDicksonModif}
 \begin{equation*}
   \forall_{f,g,n} \ex_k(
   I(n) \le k \land
   D'(k)).
 \end{equation*}
\end{prop}

\begin{proof}
 By general induction with measure function $\Phi$.  We fix $f,g$ and
 prove progressiveness of the remaining $\forall_n \ex_k$-formula.
 Therefore we can assume as induction hypothesis that for all $m$
 with $\Phi(m)<\Phi(n)$ we have
 \begin{equation*}
   \ex_k(
   I(m) \le k \land D'(k)).
 \end{equation*}
 We must show
 \begin{equation*}
   \ex_k(
   I(n) \le k \land D'(k)).
 \end{equation*}
 By Lemma~\ref{L:DescentModif} we know $D'(I(n)) \lor
 \Phi(I(n))<\Phi(n)$.  In the first case we have $D'(I(n))$
 and can take $k \defeq I(n)$.  In the second case we apply the
 induction hypothesis to $I(n)$. It provides a $k$ with $I(I(n)) \le
 k$ and $D'(k)$.  But $I(n) \le I(I(n))$ since
 $n<I(n)$.
\end{proof}

\subsection{Formalization and extraction}
\label{SS:Extract}
The formalization\footnote{Available at
 \texttt{git/minlog/examples/arith/dickson.scm}} (in
Minlog\footnote{See \url{http://www.minlog-system.de}}) of the proof
above is now routine.  The term extracted from it is
\begin{verbatim}
[f,g,n](GRecGuard nat nat)(Phi f g)n
([n0,f1][if (cDesc f g n0) (I f g n0) (f1(I f g n0))])
True
\end{verbatim}
To explain this term we rewrite it in the notation above
\begin{equation*}
 \lambda_{f,g,n}\grec \mu n G
\end{equation*}
with measure $\mu$ and step function $G$ defined by
\begin{align*}
 \mu &\defeq \Phi,
 \\
 G(n,h) &\defeq
 \begin{cases}
   I(n) &\hbox{if $\mathtt{cDesc}(n)$, i.e., $D'(I(n))$},
   \\
   h(I(n)) &\hbox{otherwise, i.e., $\Phi(I(n)) < I(n)$},
 \end{cases}
\end{align*}
where for readability we again omit the arguments $f,g$ from $\Phi, I,
\mathtt{cDesc}$.  The functions $\Phi, I$ are defined as in
\eqref{E:PsiPhiID}, and $\mathtt{cDesc}$ is the computational content
of Lemma~\ref{L:DescentModif}:
\begin{verbatim}
[f,g,n][case (cKey f g n(f(Mini g n)max g(Mini f n)))
  ((DummyL nat ysum nat) -> True)
  (Inr nn -> 
  [case nn
    ((InL nat nat)n0 -> 
     (cNatLeLtCases boole)(g(Mini g n))(g n0)True False)
    ((InR nat nat)n0 ->
     (cNatLeLtCases boole)(f(Mini f n))(f n0)True False)])]
\end{verbatim}
Here \texttt{nn} is a variable of type $\typeN + \typeN$ with $\typeN$
the type of natural numbers, and \texttt{cNatLeLtCases}:
\begin{verbatim}
(Rec nat=>nat=>alpha=>alpha=>alpha)n
([n0,x,x0][case n0 (0 -> x0) (Succ n1 -> x)])
([n0,h,n1,x,x0][case n1 (0 -> x) (Succ n2 -> h n2 x x0)])
\end{verbatim}
is the computational content of the (simple) proof of
\begin{equation*}
 \forall_{n,m}((n \le m \to P) \to (m<n \to P) \to P)
\end{equation*}
expressing case distinction w.r.t.\ $\le$ and $<$.

$\mathtt{cKey}$ is the computational content of Lemma~\ref{L:Key}:
\begin{verbatim}
[f,g,n,n0]
[case (cFPHDisjTwo([n1]f(Succ(n+n1)))([n1]g(Succ(n+n1)))n0)
  ((DummyL nat ysum nat) -> (DummyL nat ysum nat))
  (Inr nn -> 
  Inr[case nn
       ((InL nat nat)n1 -> (InL nat nat)(Succ(n+n1)))
       ((InR nat nat)n1 -> (InR nat nat)(Succ(n+n1)))])]
\end{verbatim}
which uses \texttt{cFPHDisjTwo}:
\begin{verbatim}
[f,g,n][if (cFPHDisj(n*n)
      ([n0][if (g n0<f n0)
               (f n0*f n0+g n0)
               (g n0*g n0+g n0+f n0)]))
  ([ij](DummyL nat ysum nat))
  ([n0]
   Inr[if (cCodeSqFill(f n0)(g n0)n)
          ((InL nat nat)n0)
          ((InR nat nat)n0)])]
\end{verbatim} 
which in turn depends on \texttt{cCodeSqFill}:
\begin{verbatim}
[n,n0,n1](Rec nat=>nat=>boole)n([n2]False)
([n2,(nat=>boole),n3]
  [case n3 (0 -> True) (Succ n -> (nat=>boole)n)])
n0
\end{verbatim}
and \texttt{cFPHDisj}:
\begin{verbatim}
[n](Rec nat=>(nat=>nat)=>nat@@nat ysum nat)n
([f](InR nat nat@@nat)0)
([n0,d,f]
  [let n1
    [if (f(Succ n0)<=f(Maxi f n0)) (Maxi f n0) (Succ n0)]
    [if (Succ n0<=f n1)
     ((InR nat nat@@nat)n1)
     [if (d([n2][if (n2<n1) (f n2) (f(Succ n2))]))
      ([ij]
       (InL nat@@nat nat)
       [if (right ij<n1)
         ij
         ([if (left ij<n1)
              (left ij)
              (Succ left ij)]@Succ right ij)])
      ([n2]
       [if (n2<n1)
         ((cNatLeCases nat@@nat ysum nat)(f n2)(f n1)
         ((InL nat@@nat nat)(0@0))
         ((InL nat@@nat nat)(n2@n1)))
         ((cNatLeCases nat@@nat ysum nat)(f(Succ n2))(f n1)
         ((InL nat@@nat nat)(0@0))
         ((InL nat@@nat nat)(n1@Succ n2)))])]]])
\end{verbatim}


To summarize, we have extracted a function which takes two functions
$f,g$ (suppressed for readability) and a number $n$ and yields a
bound.  Notice that already with $n=0$ we obtain the desired bound for
Dickson's lemma.  However, the inductive argument requires the general
formulation.

Our extracted bound $B(n) \defeq \grec \Phi n G$ satisfies
\begin{align*}
 B(n) = \grec \Phi n G &= G n (
 \lambda_m \ifthenels{\Phi m<\Phi n}{\grec \Phi m G}{\nullterm})
 \\
 &=\begin{cases}
   I(n) &\hbox{if $D'(I(n))$},
   \\
   B(I(n)) &\hbox{if $\Phi(I(n)) < I(n)$}.
 \end{cases}
\end{align*}
by Lemma~\ref{L:DescentModif}, which also guarantees termination:
$B(n)$ will call itself at most $I(n)$ times.  As long as the
iterations $I(n)$, $I^2(n)$, \dots, $I^m(n)$ decrease w.r.t.\ the
measure $\Phi$, the next iteration step is done.  However, as soon as
Lemma~\ref{L:DescentModif} goes to its \inquotes{left} alternative
(i.e., $D'(I(n))$ holds), $I(n)$ is returned.  Hence this extracted
bound differs from the \inquotes{guessed} one in
Proposition~\ref{P:BoundDickson} in that it does not iterate $I$ a
prescribed number of times (${f_0+g_0+1}$ many) at $0$, but stops when
allowed to do so by the outcome of Lemma~\ref{L:DescentModif}.

\bibliography{bd}
\bibliographystyle{abbrv}

\end{document}